
\input amstex
\documentstyle{amsppt}
\magnification \magstep1
\parskip10pt
\parindent.3in
\pagewidth{5.2in}
\pageheight{7.2in}
\NoRunningHeads

\centerline {\bf A Matsusaka-type Theorem for Surfaces}
\vskip.4cm
\centerline {\smc Guillermo Fern\'andez del Busto \footnote{Partially
supported by UNAM}}

\document

\def\Pj{\Bbb P}

\def\fp{\quad\quad\square}
\def\lrou{\lceil\!}
\def\rrou{\!\rceil}
\def\dsum{\displaystyle{\sum_{i=1}^r}}
\def\bar{\overline}
\def\k{\kappa}
\def\bk{\beta(k)}
\def\m{\frak m}
\def\O{\Cal O}

\vskip.8cm
\subheading{Introduction}

\vskip.4cm
The purpose of this note is to give an effective version of
Matsusaka's theorem on a smooth projective surface. Matsusaka's
general theorem states that given an ample line bundle $A$ on a smooth
projective variety (of arbitrary dimension) $X$ there exists a constant
$k_0$, depending only on the coefficients of the Hilbert polynomial of $A$,
such that $kA$ is very ample for $k\geq k_0$. In case $X$ is a curve of genus
$g$, then $A$ is ample if and only if $\deg (A)>0$, and in this case $kA$ is
very ample if $k\geq \frac{2g+1}{\deg (A)}$; for an arbitrary surface and
ample divisor $A$ (in any characteristic) the theorem was proven by
Matsusaka and Mumford in \cite{MM} and the general case by
Matsusaka in \cite{Ma}.

\vskip.4cm
The original proof of Matsusaka's theorem depends on the
boundedness of certain invariants of certain varieties and divisors in a
bounded family; the constant $k_0$ cannot be effectively computed.
Matsusaka's original approach is non-cohomological, but as pointed out by
Lieberman and Mumford~\cite{LM}, it suffices to find an integer $k_0$ such
that for any polarized variety $(X,A)$ with given Hilbert polynomial, and any
pair of points $x,y \in X$
$$H^1(X,\O _X(kA) \otimes {\m _x} \otimes {\m _y})=0
    \phantom{..}\hbox{for} \phantom{..}k\geq k_0 \phantom{.}.$$

\vskip.4cm
Recently, Siu \cite{S} gave an effective version of Matsusaka's theorem in
all dimensions. His method uses the strong Morse inequality, the numerical
criterion for very ampleness of Demailly and Nadel's vanishing theorem.
Unfortunately, his lower bound depends doubly exponentially on the dimension
of the variety; in this sense Siu's bound is more of theoretical interest than
of practical utility.  In this note, we focus only on the case of surfaces,
and obtain a result which is esentially optimal, as is shown by an
example of Xiao. Our version of Matsusaka's theorem is the following:

\noindent
\proclaim {Theorem}
Let $A$ be an ample divisor on a nonsingular projective algebraic surface $X$.
If
$$k>\frac 12 \Bigl[\frac {(A\cdot (K_X+4A)+1)^2}{A^2}+ 3\Bigr] $$
then $kA$ is very ample.
\endproclaim

\vskip.8cm
Following Siu, the argument starts by recalling a classical lemma
guaranteeing that some multiple of the difference of suitable divisors has a
lot of sections. Next,  we use a technique of Ein and
Lazarsfeld ---very much motivated by some of Demailly's analytic ideas
(c.f.~\cite{D}) ---
to produce some auxiliary divisors with almost isolated singularities.
Then, using the cohomological techniques pioneered by Kawamata, Reid and
Shokurov in connection with the minimal model program (c.f. \cite{CKM} or
\cite{KMM}), we are able to conclude with the Kawamata-Viehweg vanishing
theorem for
fractional divisors. Along the way we give a greatly simplified proof of the
main lemma of Lazarsfeld, Ein and Nakamaye \cite{L} required for this
technique, which potentially opens the door to higher dimensional extensions.

\vskip.4cm
At the end of this note, we reproduce an example of Gang Xiao,
kindly communicated
to me by Lawrence Ein, which shows that this version of Matsusaka's theorem is
esentially optimal.

\vskip.4cm
This note is part of my Ph.D. thesis at UCLA, and I would like to
thank Rob Lazarsfeld for his guidance and encouragement.
I would further like to thank Lawrence Ein for valuable discussions, Daniel
Huybrechts for valuable sugestions,
as well as Gang Xiao for allowing me to reproduce his example
in this note.

\vskip.8cm
\subheading{\S 1. Ein-Lazarsfeld-Nakamaye's lemma}

\vskip.4cm
Consider a nonsingular projective surface $X$ and let $B$ be a big divisor
on $X$, i.e. suppose there exists a constant $\beta >0$ such that
$h^0(X,\O_X(nB))\geq \frac{n^2}2 \beta$ for $n\gg 0$ (recall that if
in addition $B$ is nef, this condition is equivalent to $B^2>0$).

\vskip.4cm
For $n>0$, let
$$nB=M_n+F_n$$
be the decomposition of $nB$ in it's {\it moving part} $M_n$ and it's
{\it fixed part} $F_n$, so that $M_n$ is nef (and big) and

$$h^0(X, \O_X(M_n))=h^0(X, \O_X(nB)) \phantom{.}. \tag 1$$

Note that as a consequence of Bertini's theorem, the general element of the
complete linear system $|M_n|$ is reduced.

\vskip.4cm
Ein-Lazarsfeld-Nakamaye's technique to construct a divisor linearly equivalent
to $nB$ with (almost) isolated singularities, consists in finding first a lower
bound on the self-intersection of $M_n$ in terms of the number of sections of
$\O_X(nB)$, and then use the Hodge index theorem to find an upper bound on the
coefficients of the fractional divisor
$\frac 1n F_n$. More precisely, we first need the following (c.f.~\cite{L}):

\vskip.8cm
\proclaim{Lemma (Ein-Lazarsfeld-Nakamaye)}
Suppose that for $n\gg 0$, there is a constant $\beta >0$ with
$$h^0(X,\O_X(nB))\geq \frac {n^2}2 \beta +o(n) \phantom{.}.$$
Then $M_n^2 \geq n^2\beta +o(n)$ for $n\gg 0$.
\endproclaim

\demo{Proof}
Let $H$ be a very ample divisor on $X$, with the property that
$K_X+H$ is very ample. A general element of the linear system $|K_X+H|$
defines, for every $n$, an inclusion $\O_X(M_n)\subseteq \O_X(K_X+H+M_n)$,
so in particular
$$h^O(X,\O_X(K_X+H+M_n))\geq h^0(X,\O_X(M_n)) \phantom{.}.$$
Since $H+M_n$ is nef and big, the Kawamata-Viehweg vanishing theorem implies
that $\chi (X,\O_X(K_X+H+M_n)) = h^0(X,\O_X(K_X+H+M_n))$, and hence it
follows from Riemann-Roch that
$$\frac 12 ((K_X+H+M_n)\cdot (H+M_n))+ o(n) \geq h^0(X,\O_X(M_n))
\phantom{.}.$$
But $M_n \cdot (K_X+H)\sim o(n)$, in fact, since $K_X+H$ is very ample
$$nB\cdot (K_X+H) \geq M_n \cdot (K_X+H) \phantom{.};$$
the result now follows from~(1).$\fp$
\enddemo

\vskip.8cm
Now suppose that $N$ is an arbitrary nef and big divisor on $X$. Then
$F_n\cdot N\geq 0$ and the Hodge
index theorem implies that
$$\aligned
        \frac 1n F_n\cdot N &=B\cdot N -\frac 1n M_n\cdot N \\
        &\leq B\cdot N-\sqrt{\frac {M_n^2}{n^2} } \sqrt{N^2} \phantom{.},
  \endaligned$$
so in the situation of the previous lemma
$$ 0\leq \frac 1n F_n\cdot N
    \leq B\cdot N-\sqrt{\beta} \sqrt{N^2} \tag 2$$
and in particular an upper bound on the right hand side of~(2) imposes
numerical conditions on the coefficients of $F_n$ and in consequence on the
divisor $[\frac 1n F_n]$.

\vskip.8cm
\subheading{\S 2. Proof of theorem}

\vskip.4cm
Following Beltrametti and Sommese~\cite{BS}, recall that a divisor $L$
is {\it $\k$-jet ample} if, given any $r$ positive
integers $\k_1,\dots ,\k_r$ with
$\k+1=\dsum  \k_i$, and any
$r$ distinct points $\{ x_1, \dots ,x_r \} \subseteq X$, the evaluation map
$$H^0(X, \O_X(L))\longrightarrow H^0(X, \O_X(L) \otimes \O_Z)$$
is surjective, where $Z=\k_1 x_1 +\cdots + \k_r x_r$. In particular
$L$ is globally generated (resp. very ample) if and only if $L$
is $0$-jet ample (resp. $1$-jet ample).

\vskip.4cm
The purpose now is to prove the following generalized effective version of
Matsusaka's Big theorem.

\vskip.8cm
\proclaim{Theorem}
Let $A$ be an ample divisor on a nonsingular projective
surface $X$. If
$$k>\frac 12 \Bigl[\frac {(A\cdot (K_X+4A)+1)^2}{A^2} + (\k ^2+4\k -2)\Bigr]$$
then $kA$ is $\k$-jet ample.
\endproclaim

\vskip.4cm
For the proof, let $\k_1,\dots ,\k_r$ be $r$ positive integers with
$\k+1=\dsum \k_i$, and let
$\{ x_1, \dots ,x_r \} \subseteq X$ be any $r$ distinct points. Let $Z$ denote
the 0-cycle $Z=\dsum \k_i x_i$ and let
$$\goth m_Z =
   \goth m_{x_1}^{\k_1} \otimes \cdots \otimes \goth m_{x_r}^{\k_r}$$
where $\goth m_{x_i}$ denotes the maximal ideal at $x_i$. Then for $kA$ to
be $\k$-jet ample it suffices to show that

$$H^1(X,\O_X(kA) \otimes \goth m_Z )=0 \phantom{.}. \tag 3$$

\vskip.4cm
To this end, let $f:Y\rightarrow X$ be the blow-up of $X$ at $x_1,\dots , x_r$,
with corresponding exceptional divisors $E_1,\dots ,E_r\subseteq Y$.
Using the Leray spectral sequence, ~(3) is then equivalent to
$$H^1(Y, \O_Y(f^*(kA)-\dsum \k_i E_i ))=0 \phantom{.}. \tag 4$$

\vskip.4cm
Let $k$ be a positive integer.
In order to apply Ein-Lazarsfeld-Nakamaye's lemma, we
need to find a lower bound on the number of sections of $n(kA-K_X)$ for
$n\gg 0$. Writing
$$kA-K_X=(k+4)A-(K_X+4A) \phantom{,},$$
the question reduces to the case of the difference of two ample divisors (a
result of Fujita~\cite{F} asserts that $K_X+4A$ is ample whenever A is
ample). In the case of surfaces, this is (an easy) well known consequence of
Riemann-Roch (c.f. \cite{Mu}).

\vskip.8cm
\proclaim{Lemma 1}
Let $D$ and $E$ be ample divisors on a nonsingular projective
surface $X$. If $D^2-2D\cdot E>0$ then $h^0(X,\O_X(n(D-E))) \not = 0$ for
$n\gg 0$. In fact, for $n\gg 0$
$$h^0(X,\O_X(n(D-E)))\geq
             \frac {n^2}2 (D^2-2D\cdot E)+o(n)\phantom{.}. \tag 5$$
\endproclaim

\vskip.8cm
In particular, if $k$ is a positive integer such that
$$(k+4)A^2-2A\cdot (K_X+4A)>0 \phantom{,},\tag 6$$
then for $n\gg 0$ we have that
$$h^0(X,\O_X(n(kA-K_X)))\geq \frac {n^2}2 ((k+4)^2A^2-2(k+4)A\cdot
(K_X+4A))+o(n) \phantom{.}.$$
Write $\bk=(k+4)A^2-2A\cdot (K_X+4A)$ and
$\bar{\bk}=(k+4)\bk$.

\vskip.4cm
Next we apply Ein-Lazarsfeld-Nakamaye's lemma on the blow-up $Y$ of $X$ at
$x_1,\dots ,x_r$. For this, consider the divisor
$$B=f^*(kA-K_X)-\dsum (\k_i +1)E_i \phantom{.}.$$

\vskip.8cm
\proclaim{Lemma 2}
With notation as above, if $k$ is such that $\bk >0$ then
$$h^0(Y,\O_Y(nB))\geq
             \frac {n^2}2 (\bar{\bk}-(\k +2)^2)+0(n) \phantom{.}. \tag 7$$
\endproclaim

\demo{Proof}
Recall that the integers $\k_1,\dots ,\k_r$ are such that
$\dsum \k_i =\k +1$; under this condition, a direct computation shows
(c.f.~\cite{BS}) that
$$(\k +2)^2 \geq
           \dsum (\k_i+1)^2 +(r-1)^2 \phantom{.}. \tag 8$$
On the other hand, we have that
$$h^0(Y,\O_Y(nB)) \geq
    h^0(X,\O_X(n(kA-K_X))-\dsum \ell (\O_X / \goth m_{x_i}^{n(\k_i +1)} )$$
and that
$\ell (\O_X / \goth m_{x_i}^{\k_i} )= {{n(\k_i+1)+1}\choose{2}} $,
so it follows from~(5) and~(8) that
$$\eqalign{
        h^0(Y,\O_Y(nB))&\geq
           \frac {n^2}{2} (\bar{\bk} -\dsum (\k_i+1)^2) +o(n) \cr
                       &\geq
           \frac {n^2}{2} (\bar{\bk} -(\k+2)^2) +o(n) \phantom{.}. \fp \cr}$$
\enddemo

\vskip.4cm
In particular for $B$ to be big it suffices by~(7) that

$$\bar{\bk}-(\k+2)^2 >0 \phantom{.}. \tag 9$$
Suppose that $B$ is big, that is, suppose that $k$ is such that
(9) holds.  Then Ein-Lazarsfeld-Nakamaye's lemma implies that
$$\frac 1n F_n\cdot f^*A \leq
   A\cdot (kA-K_X)-\sqrt{\bar{\bk} -(\k+2)^2} \sqrt{A^2} \phantom{.}.$$

\vskip.8cm
\proclaim{Lemma 3}
Suppose that $k$ is such that
$$k>\frac 12 \Bigl[\frac {(A\cdot (K_X+4A)+1)^2}{A^2} +(\k ^2+4\k -4) \Bigr]
       \phantom{.}.$$
Then $\bar{\bk}-(\k+2)^2 >0$ and
$$A\cdot (kA-K_X)-\sqrt{\bar{\bk} -(\k+2)^2} \sqrt{A^2} <1\phantom{.}.$$
\endproclaim

\demo{Proof}
This is a direct computation.$\fp$
\enddemo

\vskip.8cm
Now suppose that
$$k>\frac 12 \Bigl[\frac {(A\cdot (K_X+4A)+1)^2}{A^2} +(\k ^2+4\k -4)\Bigr]
       \phantom{.}.$$
Then $\frac 1n F_n\cdot f^*A <1$ and hence the irreducible components of
the divisor $[\frac 1n F_n]$ are exceptional, that is,
$$[\frac 1n F_n]=\dsum \eta _iE_i \quad (\eta _i\geq 0) \phantom{.}.
  \tag 10$$

\vskip.4cm
Fix $n\gg 0$, and let $D$ be a general divisor on the linear system $|nB|$ so
that if $D=M+F_n$ with $F_n$ the fixed part of $nB$ then $M$ is reduced.
The divisor
$$f^*((k+1)A-K_X)-\dsum (\k_i+1)E_i-\frac 1n D $$
being numerically equivalent to $f^*A$ is nef and big. Using~(10) and
the fact that $M$ is reduced, the Kawamata-Viehweg vanishing theorem then
implies the vanishing of the higher cohomology groups of the divisor
$$K_Y+\lrou f^*((k+1)A-K_X)-\dsum (\k_i+1)E_i-\frac 1n D\rrou =
    (k+1)f^*A-\dsum (\k_i+\eta _i)E_i\phantom{,},$$
and in particular we get that
$$H^1(Y,f^*\O_X((k+1)A)\otimes
  \O_Y(-\dsum (\k_i+\eta _i)E_i))=0 \phantom{.}. \tag 11$$
But~(11) implies~(4) and it follows that $(k+1)A$ is $\k$-jet ample.
This completes the proof of the theorem.

\vskip.8cm
As an immediate consequence, we get the following effective version of
Matsusaka's theorem:

\vskip.4cm
\proclaim{Corollary} Let $A$ be an ample divisor on a nonsingular projective
surface. If
$$k>\frac 12 \Bigl[\frac {(A\cdot (K_X+4A)+1)^2}{A^2} + 3\Bigr] \tag 12$$
then $kA$ is very ample.
\endproclaim

\vskip.8cm
The following example of Xiao shows that the lower bound on (12) is
esentially optimal.

\vskip.4cm
\proclaim{Example (Xiao)}
{\rm Let $X_n$ be the ruled surface over $\Pj ^1$ defined by
$\O_{\Pj ^1}\oplus \O_{\Pj ^1}(-n)$ with section $C_0$ corresponding to
$\O_{X_n}(C_0)=\O_{X_n}(1)$ and fiber $f$. Let
$C_1\sim C_0+nf$ be a section with $C_0\cap C_1=\emptyset$ and let $D$ be an
irreducible nonsingular curve on the linear system $|2C_1+nf|$
(c.f.~\cite{Ha, III.2}). For the surface $X$ consider the double cover
$\phi \colon X\rightarrow X_n$ of $X_n$ branched along $D$, and for the ample
divisor $A$ consider $\phi ^*(C_1+f)$.
Note that, for $n\geq 6$, $\phi ^*(C_0)\subseteq X$
is a hyperelliptic curve of genus $\frac n2 -1$ and that $A|_{\phi ^*(C_0)}$
is an ample divisor of degree 2 on $\phi ^*(C_0)$; in particular, if $kA$
is very ample then
$$kA\cdot \phi ^*(C_0) > \frac n2 +1 \phantom{.}. \tag13 $$
On the other hand $A^2=2(n+2)$ and $A\cdot K_X=n-6$, which implies by~(13)
that if $kA$ is very ample then
$$k>\frac 12 \Bigl[ \frac {(A\cdot K_X)^2+8(A\cdot K_X)}{A^2} +1\Bigr]
    \phantom{.}.$$}
\endproclaim

\vskip.4cm
\proclaim{Remark}
{\rm Using the same method as above, with a different lower bound on the number
of sections of $n(kA-K_X)$ (namely $\bar{\bk} =(kA-K_X)^2$), one can show the
following effective variant of Matsusaka's theorem:}
\endproclaim

\vskip.4cm
\proclaim{Theorem$^*$}
Let $A$ be an ample divisor on a nonsingular projective surface. If
$$k>\frac 12 \Bigl[\frac {(A\cdot K_X +1)^2}{A^2} -K_X^2 +(\k ^2+4\k +6)\Bigr]
    \tag 13$$
then $kA$ is $\k$-jet ample.
\endproclaim

\vskip.8cm

\settabs\+ Institute for Advanced Study, School of Mathematics \cr
\+ School of Mathematics \cr
\+ Institute for Advanced Study \cr
\+ Princeton, NJ 08540 \cr
\+ E-mail: gfb\@math.ias.edu \cr


\vskip.8cm
\widestnumber\key{KMM}

\Refs
\ref\key BS
    \by M. Beltrametti and A. Sommese
    \pages 355--376
\endref

\ref\key CKM
    \by H. Clemens, J. Koll\'ar and S. Mori
    \book Higher dimensional Complex Geometry
    \publ Asterisque {\bf 166}
    \yr 1988
\endref

\ref\key D
    \by J.-P. Demailly
    \pages 323-374
    \paper A numerical criterion for very ample line bundles
    \jour J. Diff. Geom {\bf 37}
    \yr 1993
\endref

\ref\key F
    \by T. Fujita
    \pages 167--178
    \paper On polarized manifolds whose adjoint bundles are not semipositive
    \inbook Algebraic Geometry, Sendai
    \bookinfo Adv. Studies in Pure Math {\bf 10}
    \publ Kinokuniya-North-Holland
    \yr 1987
\endref

\ref\key Ha
    \by R. Hartshorne
    \book Algebraic Geometry
    \bookinfo Graduate Text in Mathematics {\bf 52}
    \publ Springer Verlag, New York
    \yr 1977
\endref

\ref\key KMM
    \by Y. Kawamata, K. Matsuda and K. Matsuki
    \pages 283--360
    \paper Introduction to the minimal model problem
    \inbook Algebraic Geometry, Sendai
    \bookinfo  Adv. Studies in Pure Math {\bf 10}
    \publ Kinokuniya-North-Holland
    \yr 1987
\endref

\ref\key L
    \by R. Lazarsfeld
    \paper Lectures on Linear Series
    \inbook Park City/IAS Math. Series Vol. 3
    \toappear
\endref

\ref\key LM
    \by D. Lieberman and D. Mumford
    \pages 513--530
    \paper
    \jour Proceedings of Symposia in Pure Mathematics Vol. 29
    \yr 1975
\endref

\ref\key Ma
    \by T. Matsusaka
    \pages 283--292
    \paper On canonically polarized varieties (II)
    \jour  Amer. J. Math.
    \vol 92
    \yr 1970
\endref

\ref\key MM
    \by T. Matsusaka and D. Mumford
    \pages 668--684
    \paper Two fundamental theorems on deformations of polarized varieties
    \jour Amer. J. Math
    \vol 86
    \yr 1964
\endref

\ref\key Mu
    \by D. Mumford
    \book Lectures on Curves on an Algebraic Surface
    \publ  Princeton Univ. Press
    \yr 1966
\endref

\ref\key S
    \by Y-T Siu
    \paper An Effective Matsusaka Big Theorem
    \toappear.
\endref
\endRefs
\enddocument